\documentstyle[amsfonts,12pt]{article}
\input amssym.def

\newcommand{\Dif}{{\rm Diff}}
\newcommand{\tq}{q_u}
\newcommand{\tE}{\widetilde{E}}
\newcommand{\tP}{\widetilde{P}}
\newcommand{\tstar}{\tilde{*}}
\newcommand{\EL}{End_{L_N}}

\newcommand{\EV}{End_V}
\newcommand{\HV}{End(V)}
\newcommand{\hA}{{\hat{A}}}

\newcommand{\tD}{\widetilde{D}}
\newcommand{\tA}{\widetilde{A}}
\newcommand{\tn}{\widetilde{\nabla}}
\newcommand{\tpa}{\widetilde{\partial}}
\newcommand{\nr}{\widetilde{r}}

\newcommand{\tV}{\widetilde{V}}
\newcommand{\tla}{\widetilde{\lambda}}

\newcommand{\al}{{\alpha}}
\newcommand{\la}{{\lambda}}
\newcommand{\h}{{\hbar}}
\newcommand{\bul}{{\bullet}}

\newcommand{\mg}{{\mathfrak{g}}}
\newcommand{\mh}{{\mathfrak{h}}}

\newcommand{\mgl}{{\mathfrak{gl}}}
\newcommand{\msp}{{\mathfrak{sp}}}

\newcommand{\ph}{{\varphi}}
\newcommand{\om}{{\omega}}
\newcommand{\bom}{\overline{\omega}}
\newcommand{\Om}{{\Omega}}

\newcommand{\si}{{\sigma}}

\newcommand{\ve}{{\varepsilon}}

\newcommand{\G}{{\Gamma}}

\newcommand{\pa}{{\partial}}

\newcommand{\M}{{\cal M}}

\newcommand{\cD}{{\cal D}}
\newcommand{\cA}{{\cal A}}

\newcommand{\cW}{{\cal W}}

\newcommand{\cE}{{\cal E}}

\newcommand{\cO}{{\cal O}}
\newcommand{\cV}{{\cal V}}

\newcommand{\bbC}{{\Bbb C}}
\newcommand{\bbR}{{\Bbb R}}

\newcommand{\bbA}{{\Bbb A}}
\newcommand{\La}{{\Lambda}}

\newcommand{\n}{{\nabla}}

\newcommand{\Te}{\Theta}

\newcommand{\de}{{\delta}}
\newcommand{\D}{{\Delta}}

\date{}
\newtheorem{defi}{Definition}
\newtheorem{lem}{Lemma}
\newtheorem{teo}{Theorem}
\newtheorem{cor}{Corollary}
\newtheorem{pred}{Proposition}
\begin{document}

\vspace{-10cm}
\begin{flushright}
 \begin{minipage}{1.2in}
 ITEP-TH-29/04
 \end{minipage}
\end{flushright}
\vspace{1.1cm}

\begin{center}
{\Large\bf A SIMPLE  ALGEBRAIC PROOF OF\\[0.5cm]
 THE ALGEBRAIC INDEX THEOREM.}\\[1cm]
PoNing Chen and
Vasiliy Dolgushev\footnote{On leave of
absence from: ITEP (Moscow)}\\[0.5cm]
{\it Department of Mathematics, MIT,} \\
{\it 77 Massachusetts Avenue,} \\
{\it Cambridge, MA, USA 02139-4307,}\\
E-mails: pnchen@MIT.EDU, vald@MIT.EDU\\[1cm]
\end{center}

\begin{abstract}
In math.QA/0311303 B. Feigin, G. Felder, and B. Shoikhet
proposed an explicit formula for the trace density
map from the quantum algebra of functions on an
arbitrary symplectic manifold $\M$ to the top degree
cohomology of $\M$. They also evaluated this map
on the trivial element of $K$-theory of the
algebra of quantum functions. In our paper we
evaluate the map on an arbitrary element
of $K$-theory, and show that the result is
expressed in terms of the $\hA$-genus
of $\M$, the Deligne-Fedosov class of the quantum algebra,
and the Chern character of the principal
symbol of the element. For a smooth (real)
symplectic manifold (without a boundary),
this result implies the Fedosov-Nest-Tsygan
algebraic index theorem.
\end{abstract}
~\\[0.3cm]
MSC-class: 19A49, 19K56.

\section{Introduction}
The famous Atiyah-Singer index theorem \cite{AS} relates the
index of an elliptic pseudo-differential operator on a compact
manifold $X$ to the Todd class of $X$ and the Chern character
of the bundle naturally associated with the symbol of
the pseudo-differential operator. It is well-known that
the algebra of pseudo-differential operators can be viewed as
a ``quantization'' of the cotangent bundle $T^{*}X$.
Starting from this point B. Fedosov \cite{F1} proposed a natural
analogue of the Atiyah-Singer index theorem for
the deformation quantization of an arbitrary symplectic
manifold. This theorem expresses the index of a compactly
supported element in $K$-theory of the quantum algebra
of functions on a symplectic manifold $\M$ via the $\hA$-genus
of $\M$, the Deligne-Fedosov class of the quantum algebra,
and the Chern character of the principal symbol of
the $K$-theory element. While the original index theorem
of Atiyah and Singer relates two integers, the Fedosov
index theorem relates two formal power series with complex
coefficients.

In paper \cite{NT} R. Nest and B. Tsygan proposed a completely
different proof of the algebraic index theorem based on
the relations between Hochschild (co)homology, cyclic
(co)homology, and Lie algebra (co)homology. In their next
paper \cite{NT1} they also showed that the analytic index theorem
for a compact manifold (without a boundary)
could be derived from the algebraic one.
In order to generalize this result to manifolds
with a boundary (or more generally, with corners) one
is in great need of {\it a local version of the algebraic
index theorem} which relates De Rham cohomology
classes, rather than formal power series of
complex numbers.

The aim of this paper is to give such a version (see theorem
\ref{index} in section $3$) using the trace density
map proposed in paper \cite{FFS} by B. Feigin, G. Felder,
and B. Shoikhet. Namely,
we show that the Feigin-Felder-Shoikhet (FFS) trace
density map sends an arbitrary element of $K$-theory of the
algebra of quantum functions on $\M$ to the top
component of the cup product of the $\hA$-genus
of $\M$, the exponent of the Deligne-Fedosov class
of the quantum algebra, and the Chern character of the principal
symbol of the element. Thus, we generalize the result
of \cite{FFS}, in which the trace density map was
evaluated at the trivial element of $K$-theory.
For the case of the cotangent bundle this computation of
\cite{FFS} was performed in \cite{CFS}.

We have to mention that
in \cite{BNT1} another version of the trace density
map was proposed. This map sends a cohomology class of the
sheaf of the periodic cyclic complex $CC^{per}_{\bul}(\cO^{quant}_X)$
of the structure sheaf $\cO^{quant}_X$ of quantum functions on
a holomorphic symplectic manifold $X$ to the cohomology
of the constant sheaf on $X$. In \cite{BNT1} it was proven that
the evaluation of this map on an arbitrary cohomology class of the
sheaf of periodic cyclic complex $CC^{per}_{\bul}(\cO^{quant}_X)$
can be expressed in terms of the principal
symbol of the class, $\hA$-genus of $X$ and the
Deligne-Fedosov class of the quantum deformation.
In \cite{BNT1} the authors refer to this result
as a Riemann-Roch theorem for deformation
quantizations. It is not hard to see that the local
version of the algebraic index theorem proven in our paper
follows from the results of P. Bressler, R. Nest and
B. Tsygan \cite{BNT1}. However, our proof is
more explicit and straightforward. For this reason we expect
that our technique is more powerful for exploring the parallel
analytic results.

We would like to mention an alternative deformation
quantization procedure due to L. Boutet de Monvel
and V. Guillemin \cite{BM-G}. This procedure allows us
to quantize an arbitrary compact symplectic manifold with
an integral symplectic form. Recently R. Melrose proposed
in \cite{RM} a beautiful idea on how one can relax the
integrality condition in the quantization of L. Boutet de Monvel
and V. Guillemin. R. Melrose suggests that using this procedure one
could derive the algebraic index theorem from the
analytic one.

The organization of the paper is as follows. In section two, we give
a reminder of Fedosov's construction \cite{F1} of the deformation quantization
with twisted coefficients on a symplectic manifold. In this section
we also recall the necessary results of Feigin, Felder, and
Shoikhet about their trace density map
proposed in \cite{FFS}.
In section three, we prove the local version
of the algebraic index theorem (see theorem {\ref{index}), using
the FFS trace density map.
In the concluding section we make a few remarks about the local
algebraic index theorem for symplectic Lie algebroids and
the versions of the local Riemann-Roch-Hirzebruch theorem.

Throughout the paper we assume the summation over repeated indices.
We assume that $\M$ is either a smooth real symplectic ma\-ni\-fold
of dimension $2n$ or a smooth affine algebraic variety (over $\bbC$)
of the complex dimension $2n$, endowed with an algebraic
symplectic form. $C(\M)$ stands, respectively, for the
algebra of smooth function or for the algebra of regular
functions on $\M$. The notation $S_k$ is reserved  for
the group of permutations of $k$ elements.
We omit the symbol $\wedge$ referring to a local basis
of exterior forms, as if we thought of $dx^i$'s as anti-commuting variables.
For a vector space $\mh$ we denote by $S^j\mh$ the subspace
of monomials of degree $j$ in the symmetric algebra $S\mh$ of $\mh$\,.
Finally, we always assume that a nilpotent linear operator is
the one whose second power is vanishing.

{\bf Acknowledgments.} This work was done while the first author
was a student and the second author was a mentor in
the Summer Program for Undergraduate Research.
We both express sincere thanks to Hartley Rogers who
organizes this remarkable program at M.I.T.
We thank Richard Melrose and all participants
of his Spring analysis seminar for their interest and
stimulating discussions. The second author also thanks
Alberto Cattaneo and Giovanni Felder for useful
conversations. We are grateful to Teal Guidici for criticisms
concerning our English.  This work is partially
supported by the NSF grant DMS-9988796, the Grant
for Support of Scientific Schools NSh-1999.2003.2,
the grant INTAS 00-561 and the grant CRDF
RM1-2545-MO-03.

\section{Preliminaries}
In this section we review the Fedosov
deformation quantization of endomorphisms
of a vector bundle $V$ over a symplectic
manifold $\M$ (See \cite{F1}, section $5.3$).
We also recall the Felder-Feigin-Shoikhet
construction \cite{FFS} of the trace density map.

\subsection{The Feigin-Felder-Shoikhet cocycles}
First, we recall that
for any associative algebra $\cA$ with a unit
over a field $K$ (of characteristic zero) we have a chain map
$\phi^N$ from the Hochschild cochain complex $C^{\bullet}(\cA,\cA^*)$
with coefficients in the dual module $\cA^*$ to
the Lie algebra cochain complex $C^{\bullet}(\mgl_N (\cA);
\mgl_N(\cA)^{*})$ with coefficients in $\mgl_N(\cA)^*$
$$
\phi^N(\psi)(M_1\otimes a_1 ,\dots,M_k\otimes a_k)(M_0\otimes a_0)
$$
\begin{equation}\label{phi}
=\frac1{k!}\sum_{\nu\in S_k}(-)^{\nu}
\psi(a_{\nu(1)}\otimes \dots \otimes a_{\nu(k)})(a_0)\,tr(M_0 M_{\nu(1)}\dots
M_{\nu(k)})\,,
\end{equation}
where $M_l\in \mgl_N(K)$, $a_l\in \cA$, and
$\psi\in C^k(\cA, \cA^*)$\,.

Let $\cV$ be a $2n$-dimensional vector space (over $\bbC$)
endowed with a symplectic form $B$.
Let $\{y^1, \dots, y^{2n}\}$ be a basis
in $\cV$. We denote by $||B^{ij}||$,
$(i,j=1, \dots , 2n)$ the matrix
$$
B^{ij} = B(y^i, y^j)
$$
of the form $B$ in this basis.
Let $\ve$ be the completely antisymmetric Liouville
tensor whose components $\ve^{i_1\dots i_{2n}}$
are defined in the basis $\{ y^1, \dots, y^{2n}\}$
as follows
\begin{equation}
\label{ve}
\ve^{i_1 \dots i_{2n}} =
\frac{(-1)^n}{n!} \sum_{\nu\in S_{2n}}
(-)^{\nu} B^{i_{\nu(1)} i_{\nu(2)}} \dots
B^{i_{\nu(2n-1)} i_{\nu(2n)}}\,.
\end{equation}

\begin{defi}\label{fWa}
The Weyl algebra $\cW$ associated with the symplectic
vector space $\cV$ is the vector space $\bbC[[\cV]]((\h))$
of the formally completed symmetric algebra of $\cV$
equipped with the following
(associative) multiplication
\begin{equation}
\label{circ}
(a \circ b)(y,\h)=
\exp\left(\frac{\h}{2}B^{ij}\frac{\pa}{\pa y^i}\frac{\pa}{\pa z^j}
\right)a(y,\h)b(z,\h)|_{y=z}.
\end{equation}
\end{defi}
One can easily see that the multiplication
defined by (\ref{circ}) does not depend on the
choice of a basis in $\cV$. We view $\cW$ as an algebra
over the field $\bbC((\h))$\,.

The Weyl algebra $\cW$ is naturally filtered with respect
to the degree of monomials $2[\h]+[y]$ where
$[\h]$ is a degree in $\h$ and $[y]$ is a degree in $y$
\begin{equation}
\label{filtr}
\begin{array}{c}
\displaystyle
\dots \subset \cW^1  \subset \cW^0  \subset \cW^{-1}
\dots \subset \cW\,, \\[0.3cm]
\displaystyle
\cW^m = \{ a = \sum_{2k+p\ge m}
\h^k a_{k; i_1 \dots i_p} y^{i_1} \dots y^{i_p}\}
\end{array}
\end{equation}
This filtration defines the $2[\h]+[y]$-adic topology in $\cW$.

In \cite{FFS} Feigin, Felder, and Shoikhet proposed an
explicit expression for the $(2n)$th Hoch\-schild (continuous)
cocycle of the Weyl algebra $\cW$ with coefficients
in the dual module $\cW^*$. Their formula is reminiscent
of Kontsevich's construction of the structure maps
of his celebrated formality quasi-isomorphism \cite{K}.
However, unlike the integrals in
Kontsevich's construction \cite{K} the coefficients
entering the formula of Feigin, Felder, and Shoikhet
are rational.

Following \cite{FFS} we start with the $2n$-th simplex
$\Delta_{2n}=\{(u_1,\dots, u_{2n})\in\bbR^{2n}\, |\, 0\le u_1\dots \le u_{2n} \le 1\}$
with the standard orientation and denote by $B_{bc}$ the action
of the form $B$ on the $b$-th and $c$-th
components of $\cW^{\otimes (2n+1)}$
$$
B_{bc}(a_0\otimes\dots \otimes a_{2n})=
$$
$$
B^{ij} (a_0\otimes\dots\otimes\frac{\pa a_{b}}{\pa y^{i}}\otimes\dots
\otimes\frac{\pa a_c}{\pa y^{j}}\otimes\dots\otimes a_{2n})\,,
\qquad a_b \in \cW\,.
$$
Furthermore, we denote by $\pi_{2n}$
the action of the Liouville tensor (\ref{ve})
$$
\pi_{2n}(a_0\otimes\dots \otimes a_{2n})=
\ve^{i_1 \dots i_{2n}}
a_0\otimes\frac{\pa a_1}{\pa y^{i_1}}\otimes\dots
\otimes\frac{\pa a_{2n}}{\pa y^{i_{2n}}}.
$$
Finally, if we denote by $\mu$ the natural projection
from $(\cW^{\otimes (2n+1)})$ onto $\bbC((\h))$
$$
\mu(a_0\otimes\dots\otimes a_{2n})=a_0(0)\dots a_{2n}(0)
$$
the Feigin-Felder-Shoikhet formula
of the $2n$-th Hochschild cocycle
$\tau_{2n} \in C^{2n}(\cW, (\cW)^*)$
can be written as \cite{FFS}
\begin{equation}
\label{tau}
\tau_{2n}(\al)(a_0)=
\mu_{2n}\Big( \int_{\Delta_{2n}}\prod_{0\le b\le c \le 2n}
e^{\frac{\h}{2}(2u_b-2u_c+1)B_{bc}}
\pi_{2n}(a_0\otimes \al)du_1\wedge\dots\wedge du_{2n}\Big)\,,
\end{equation}
where $a_0\in \cW$, $\al \in \cW^{\otimes (2n)}$, and
$u_0=0$\,.
It is easy to see that the cocycle $\tau_{2n}$ does not
depend on the choice of the basis $\{y^1, \dots, y^{2n}\}$ in
$\cV$\,.

Applying~ the~ map~ (\ref{phi})~ to~ (\ref{tau})~
we~ get~ the~ $2n$-th~ cocycle~ in~ the~ chain~ complex~
$C^{\bul}(\mgl_N(\cW), \mgl_N(\cW)^*)$ of the Lie
algebra $\mgl_N(\cW)$ with values in the dual
module $(\mgl_N(\cW))^*$
\begin{equation}
\label{Te}
\Te^N_{2n}=
\phi^N(\tau_{2n})\,:\, \wedge^{2n} (\mgl_N(\cW)) \otimes
\mgl_N(\cW) \mapsto \bbC((\h))\,.
\end{equation}
In \cite{FFS} it was shown that the latter
cocycle satisfies remarkable properties which
allow us to construct a trace density map for
a quantum algebra of functions on any symplectic manifold.
Before talking about these properties
we recall a necessary construction of
the Chern-Weil theory.

\subsection{The Chern-Weil homomorphism}
Let $\mg$ be a Lie algebra and $\mh\subset\mg$
be a subalgebra of $\mg$. Suppose that there is
an $\mh$-equivariant projection
$pr\,:\,\mg \mapsto \mh$, that is a map commuting with
the adjoint action of $\mh$ and satisfying the property
$\displaystyle pr\, |_{\mh} = Id_{\mh}\,.$
The amount by which $pr$ fails to be a Lie algebra homomorphism
is measured by the ``curvature'' $C\in Hom(\wedge^2 \mg, \mh)$
\begin{equation}
\label{curvature}
C(v,w)= [pr(v), pr(w)]- pr ([v,w])\,,
\qquad v,w \in \mg\,.
\end{equation}

Given the ``curvature'' (\ref{curvature}) of the
$\mh$-equivariant projection $pr$,
it is not hard to see that for any adjoint invariant
form $Q\in ((S^j\mh)^*)^{\mh}$ the formula
\begin{equation}
\label{Ch-W}
\chi(Q)(v_1, \dots, v_{2j}) =
\frac1{(2j)!} \sum_{\nu\in S_{2j}}
(-)^{\nu} Q(C(v_{\nu(1)}, v_{\nu(2)}), \dots,
C(v_{\nu(2j-1)}, v_{\nu(2j)}))
\end{equation}
defines a relative Lie algebra cocycle
$\chi(Q)\in C^{2j}(\mg,\mh)$. The standard argument
of the theory of characteristic classes shows that
the cohomology class $[\chi(Q)]$ of the cocycle
$\chi(Q)$ does not depend on the choice of
the projection $pr$. Thus (\ref{Ch-W}) induces
a map from the vector space $((S^j\mh)^*)^{\mh}$ to
$H^{2j}(\mg, \mh)$.
The latter map is called the {\it Chern-Weil homomorphism}.

In our case we set
$\mg=\mgl_N(\cW)$, $\mh=\mgl_N\oplus \msp_{2n}$ and define
the projection $pr$ from $\mg$ to $\mh$ by the formula
\begin{equation}
\label{pr}
pr(v)= pr_0(v) + pr_2 (v)\,,
\end{equation}
$$pr_0(v)= v\Big|_{y=0}\,, \qquad
pr_2(v)= \frac1{N}\si_2(tr (v)) I_{N}\,,$$
where $\si_2$ denotes the projection onto
the monomials of the second degree in $y$'s, $tr$ is
the ordinary matrix trace, $I_N$ is the
identity matrix of size $N\times N$, and
the Lie algebra $\msp_{2n}$ is realized as a
subalgebra of scalar matrices in
$\mgl_N(\cW)$ with values in quadratic
monomials in $\cW$\,.

Due to \cite{FFS} we have the following
\begin{teo}[Feigin-Felder-Shoikhet \cite{FFS}]
The Lie algebra cocycle $\Te^N_{2n}\in
C^{2n}(\mg,\mg^*)$ is relative with respect to the
subalgebra $\mh=\mgl_N\oplus \msp_{2n}$.
The evaluation of (\ref{Te}) on the identity matrix
$I_N\in \mg$ gives a relative Lie algebra cocycle
$$\ph\in C^{2n}(\mg, \mh)$$
$$
\ph = \Te^N_{2n}(\cdot ,\dots, \cdot, I_N)\,:\,
\wedge^{2n}(\mg) \mapsto \bbC((\h))
$$
whose cohomology class
\begin{equation}
\label{FFS}
[\ph\,]=[\chi (Q_n)]\,,
\end{equation}
coincides with the image
of the $n$-th component $Q_n\in ((S^n\mh)^*)^{\mh}$ of
the adjoint invariant form
$Q\in ((S\mh)^*)^{\mh}$
\begin{equation}
\label{Qn}
Q(X,\dots, X)= \det\Big( \frac{X_1/2\h}{sinh\,(X_1/2\h)}
\Big)^{1/2} tr\, \exp\left(\frac{X_2}{\h}\right)\,,
\end{equation}
$$
X= X_1\oplus X_2 \in \msp_{2n} \oplus \mgl_N
$$
under the Chern-Weil homomorphism (\ref{Ch-W}).
\end{teo}

\subsection{Fedosov deformation quantization with twisted
coefficients.}
Let $\M$ be a symplectic manifold of dimension $2n$.
We denote by $\om=\om_{ij}(x)dx^i dx^j$
the corresponding symplectic form and by
$$\bom = \om^{ij}(x)\frac{\pa}{\pa x^i}\wedge \frac{\pa}{\pa x^j}$$
the corresponding Poisson tensor. Here $x^i$ denote local
coordinates and the indices $i,j$ run from $1$ from $2n$\,.
For a vector bundle $V$ of rank $N$ over $\M$ we will denote
by $\HV$ the bundle of endomorphisms of $V$ and
by $End_V= \G(\M, \HV)$ the algebra of global
sections of $\HV$\,.

If the vector bundle $V$ is
endowed with a connection $\pa^V$ then
\begin{defi}
By quantization of the algebra $\EV$ we mean a construction of
an associative $\bbC((h))$-linear product in $\EV((\h))$
given by the formal power series
\begin{equation}
\label{star}
a *b = a b + \sum_{k\ge 1} \h^k B_k(a,b)\,,
\qquad a,\, b \in \EV
\end{equation}
of bidifferential operators
$B_k\,:\, \EV\otimes \EV \mapsto \EV$
such that
$$
B_1(a,b) - B_1(b,a) = \om^{ij}(x) \pa^V_i(a) \pa^V_j(b)\,.
$$
\end{defi}

In \cite{F1} B. Fedosov proposed a simple procedure
for the deformation quantization of the algebra of
endomorphisms of a vector bundle $V$ over a
symplectic manifold $\M$. The main ingredient of the construction is
the {\it Weyl algebra bundle} $W(\EV)$
whose sections are the following formal power series
\begin{equation}
a=a(x,y,h)=\sum_{k,l}\h^k a_{k; i_1 i_2\dots i_l}(x)y^{i_1}\dots
y^{i_l}\,,
\end{equation}
where $y=(y^1\dots y^{2n})$ are fiber coordinates of the tangent
bundle $T\M$, $a_{k; i_1 i_2 \dots i_l}$ represent sections of
$\HV\otimes S^l(T^{*}M)$ and the summation over $k$
is bounded below.

Multiplication of two sections of $W(\EV)$ is given by
the Weyl formula
\begin{equation}
\label{circ1}
a \circ b(x,y,\h)=\exp\left(\frac{\h}{2}\,
\om^{ij}\frac{\pa}{\pa y^i}\frac{\pa}{\pa z^j}\right)a(x,y,\h)b(x,z,\h)|_{y=z}.
\end{equation}
Notice that in the right hand side of (\ref{circ1}) the sections of
$W(\EV)$ are multiplied via the product induced
from the algebra $\EV$\,.

For any point $p\in \M$ the fiber $W_p(\EV)$ of the Weyl algebra bundle
at $p$ is isomorphic to the algebra $\mgl_N(\cW)$
of $N \times N$-matrices of the Weyl algebra $\cW$ associated
with the cotangent space $T^*_p(\M)$ at the point $p$ with the
symplectic form $\bom_p$\,.
The transition functions are realized by the adjoint
action of the group $GL_N$\,.

The filtration (\ref{filtr}) of the Weyl algebra
gives us a natural filtration of the bundle $W(\EV)$
\begin{equation}
\label{filtr-bun}
\begin{array}{c}
\displaystyle
\dots \subset W^1(\EV)  \subset W^0(\EV)  \subset
W^{-1}(\EV) \dots \subset W(\EV)\,, \\[0.3cm]
\displaystyle
\G(W^m(\EV)) = \{ a = \sum_{2k+p\ge m}
\h^k a_{k;i_1 \dots i_p}(x)y^{i_1} \dots y^{i_p}\}\,.
\end{array}
\end{equation}
This filtration defines the $2[\h]+[y]$-adic topology in the
algebra $\G(W(\EV))$ of sections of $W(\EV)$\,.

The vector space $\Om^{\bul}(W(\EV))$ of smooth exterior
forms with values in $W(\EV)$ is naturally a graded
associative algebra with the product induced by (\ref{circ1}) and the
following graded commutator
$$
[a,b]= a \circ b - (-)^{q_a q_b} b\circ a\,,
$$
where $q_a$ and $q_b$ are exterior degrees of
$a$ and $b$\,, respectively.
The filtration of $W(\EV)$ (\ref{filtr-bun}) gives us a filtration
of the algebra $\Om^{\bul}(W(\EV))$
$$
\begin{array}{c}
\dots \subset \Om^{\bul}(W^1(\EV))
\subset\Om^{\bul}(W^0(\EV))\subset  \phantom{aaaaaaa}  \\[0.3cm]
 \phantom{aaaaaaa} \subset\Om^{\bul}(W^{-1}(\EV)) \subset
\dots \subset \Om^{\bul}(W(\EV))\,.
\end{array}
$$
In what follows we refer to the algebra $\Om^{\bul}(\M)$
of exterior forms on $\M$ as an algebra embedded into
$\Om^{\bul}(W(\EV))$ via the natural map
$
\iota\,:\,\Om^{\bul}(\M) \mapsto \Om^{\bul}(W(\EV))\,,
$
which sends an exterior form $\eta\in \Om^{\bul}(\M)$
to the scalar matrix $\eta I_{N}\in \Om^{\bul}(W(\EV))$\,.

Let $\pa^s$ be a torsion free connection on $T\M$ compatible with
the symplectic structure $\om$. Using the connection $\pa^s$ and
the connection $\pa^V$ on $V$ we define the following linear
operator
$$
\n \,:\,\Om^{\bul}(W(\EV)) \mapsto
\Om^{\bul+1}(W(\EV))\,,
$$
\begin{equation}
\label{nabla}
\n=dx^i\frac{\pa}{\pa x^i}-
dx^i\G^{k}_{ij}(x)y^j\frac {\pa}{\pa y^k}+[\G_V,\,\cdot\,]\,,
\end{equation}
where $\G^k_{ij}(x)$ are Christoffel symbols of $\pa^s$
and $\G^V$ is the connection form of $\pa^V$\,.

Thanks to the compatibility of $\pa^s$ with the
symplectic structure $\om$ the operator (\ref{nabla})
is a derivation of the graded algebra
$\Om^{\bul}(W(\EV))$\,. Furthermore, a simple
computation shows that
$$
\n^2 a = \frac1{2}[R+R^V,a]\,, \qquad
\forall~ a\in \Om(W(\EV))\,,
$$
where $R^V\in \Om^2(\EV)$ is the curvature form of $\pa^V$,
$$
R= \frac1{2\h} \om_{km}(R_{ij})^m_l(x) y^k y^l dx^i dx^j,
$$
and $(R_{ij})^m_l(x)$ is the Riemann curvature
tensor of $\pa^s$\,.

\begin{defi} The~ Fedosov~ connection~ is~ a~ nilpotent
derivation of the graded algebra $\Om^{\bul}(W(\EV))$
of the following form
\begin{equation}
\label{F-conn}
D= \n + \frac1{\h}[A, \cdot\,]\,,
\qquad
A= - dx^i \om_{ij}(x)y^j + r\,,
\end{equation}
where $r$ is an element in
$\Om^1(W^2(\EV))$
\end{defi}
The flatness of $D$ is equivalent to the fact
that the Fedosov-Weyl curvature
\begin{equation}
\label{Weyl-c}
C^W= \h(R+ R^V) + 2 \n A + \frac1{\h} [A, A]
\end{equation}
of $D$ belongs to the subspace
$\Om^2(\M)((\h))\subset \Om^2(W(\EV))$\,.
A simple analysis of degrees in $\h$ and $y$
shows that $C^W$ is of the form
\begin{equation}
\label{nado}
C^W = -\om + \Om_{\h}\,,
\qquad \Om_h \in \h \Om^2(M)[[\h]]
\end{equation}
whereas the Bianchi identity $D (C^W)=0$ implies that
$\Om_h$ is a series of two-forms closed with
respect to the De Rham differential.

One can observe that the definition of the Fedosov
connection depends on the choice of the symplectic
connection $\pa^s$ and the connection $\pa^V$ on $V$.
The following proposition shows how this problem
can be remedied
\begin{pred}
\label{absor}
If $\n$ and $\tn$ are two operator (\ref{nabla}) corresponding
to the symplectic connections $\pa^s$, $\tpa^s$ and
the connections $\pa^V$, $\tpa^V$ on $V$, respectively,
then the difference $\tD-D$ between
two Fedosov connections
\begin{equation}
\label{tD}
\tD= \tn + \frac1{\h} [- dx^i \om_{ij}(x)y^j + \nr, \cdot\,]\,,
\qquad
D= \n + \frac1{\h} [- dx^i \om_{ij}(x)y^j + r, \cdot\,]\,,
\end{equation}
takes the form of the commutator
$$
\tD - D = \frac1{\h} [\D r, \cdot\,]\,,
$$
where $\D r\in \Om^1(W^2(\EV))$\,.
In other words any deviations of the connections
$\pa^s$ and $\pa^V$ can be absorbed into the
form $r\in \Om^1(W^2(\EV))$\,.
\end{pred}
{\bf Proof.} See page $151$ in \cite{F1}\,. $\Box$

Let us consider the affine subspace
$I_N\oplus \G(W^1(\EV))$ in $\G(W(\EV))$ consisting
of the sums $U=I_N + U_1\,,$ where $I_N$ is the
identity endomorphism of $V$ and $U_1$ is
an arbitrary element in $\G(W^1(\EV))$.
It is straightforward that
\begin{pred}
The affine subspace
$I_N\oplus \G(W^1(\EV)) \subset \G(W(\EV))$ is a subgroup
in the group of invertible elements of $\G(W(\EV))$\,. $\Box$
\end{pred}
We are now ready to give the following
\begin{defi}
\label{equiv}
Two Fedosov connections
$$
D = \n + \frac1{\h}[A, \cdot\,]\,, \qquad
\tD = \n + \frac1{\h}[\tA, \cdot\,]
$$
are called equivalent if there exists an element
$U\in I_N\oplus \G(W^1(\EV))$ such that
\begin{equation}
\label{equiv1}
\tD= D + [U^{-1}\circ DU,\cdot\,]
\end{equation}
or equivalently
$\tA = U^{-1}\circ A \circ U + \h U^{-1}\circ \n U\,.$
\end{defi}

Let us remark that the Fedosov connection (\ref{F-conn})
can be rewritten as
\begin{equation}
\label{DDD}
D = \n -\de + \frac1{\h}[r, \cdot\,]\,,
\end{equation}
where
\begin{equation}
\label{de}
\de= \frac1{\h}[dx^i\om_{ij}(x)y^j, \cdot\,]=
dx^i \frac{\pa}{\pa y^i}
\end{equation}
is the Koszul derivation of the algebra
$\Om^{\bul}(W(\EV))$\,.

For our purposes we will need
the homotopy operator for the Koszul
differential $\de$
\begin{equation}
\delta^{-1}a=y^k i\left( \frac \partial {\partial x^k}\right)
\int\limits_0^1 a(x,\h, ty,tdx)\frac{dt}t,  \label{del-1}
\end{equation}
where $i(\partial /\partial x^k)$ denotes the contraction
of an exterior form with the vector field
$ \partial /\partial x^k$\,, and $\delta ^{-1}$ is extended to
$\G(W(\EV))$ by zero.

Simple calculations show that
$\de^{-1}$ is indeed the homotopy operator for $\de$,
namely
\begin{equation}
a=\sigma(a) +\delta \delta ^{-1}a + \delta ^{-1}\delta a\,,
\qquad \forall~a\in\Om(W(\EV))
 \label{Hodge}
\end{equation}
where $\si$ is the natural projection
\begin{equation}
\label{sigma}
\si (a)= a \Big |_{y=0,~dx=0}\,, \qquad
a\in \Om^{\bul}(W(\EV))
\end{equation}
from $\Om^{\bul}(W(\EV))$ onto the
algebra of endomorphisms $\EV((\h))$\,.

The proof of the following theorem is contained in
section $5.3$ of \cite{F1}. (More precisely, see theorem
$5.3.3$ and remarks at the end of section $5.3$)
\begin{teo}[Fedosov, \cite{F1}]
\label{F}
If $\pa^s$ is a symplectic connection on $\M$,
$\pa^V$ is a connection on $V$ and $\Om_{\h}$
is a series of closed two-forms in $\h\Om^2(M)[[\h]]$ then
\begin{enumerate}
\item One can construct a nilpotent
derivation $D=\n + \h^{-1}[r-dx^i\om_{ij}(x)y^j, \cdot\,]$
whose Fedosov-Weyl curvature (\ref{Weyl-c}) is equal to
$C^W= -\om + \Om_{\h}$
and the element $r\in \Om^1(W^2(\EV))$ satisfies
the normalization condition
\begin{equation}
\label{norm}
\de^{-1} r = 0\,.
\end{equation}

\item Given a Fedosov connection (\ref{F-conn})
one can construct a vector space isomorphism $\la$
\begin{equation}
\label{la}
\la\,:\, \EV((\h))~ \widetilde{\rightarrow}~ \G_D(W(\EV))
\end{equation}
from $\EV((\h))$
to the algebra $\G_D(W(\EV))$
of flat sections of $W(\EV)$ with respect
to $D$. The product in $\EV((\h))$ induced via
the isomorphism $\la$ is the
desired star-product (\ref{star})

\item Two Fedosov connections $D$ and $\tD$ whose
Fedosov-Weyl curvatures represent the same De Rham
cohomology class in $H^2_{DR}(\M)[[\h]]$ are equivalent
in the sense of definition \ref{equiv}.
If the equivalence between $D$ and $\tD$ is
established by an element $U \in I_N \oplus
\G(W^1(\EV))$ then the map
$$
a \mapsto U \circ a \circ U^{-1}
$$
gives an isomorphism from the algebra $\G_{\tD}(W(\EV))$
of flat sections of $\tD$ to the algebra $\G_{D}(W(\EV))$
of flat sections of $D$\,.
\end{enumerate}
\end{teo}
{\bf Remark 1.} The construction of the Fedosov
connection is functorial in the sense that if $u$
is an isomorphism from $V$ to $\tV$ then
$$
\tD= D + [u D(u^{-1}),\cdot\,]
$$
is a Fedosov connection in $\Om^{\bul}(W(End_{\tV}))$\,.
Moreover, the corresponding isomorphisms
$$
\la\,:\, \EV((\h)) ~\widetilde{\rightarrow}~
\G_{D}(W(\EV))\,, \qquad
\tla\,:\,End_{\tV}((\h)) ~\widetilde{\rightarrow}~
\G_{\tD}(W(End_{\tV}))
$$
are related by the formula
$$
\tla(a) = u \la(u^{-1} a u) u^{-1}\,, \qquad a\in  End_{\tV}\,.
$$
{\bf Remark 2.} If $V$ is a trivial vector bundle $L_1$ of rank
$N=1$ then the isomorphism $\la$ (\ref{la}) gives a star-product $*$ in the
algebra of functions $C(\M)((\h))$\,. Due to the
result of P. Xu \cite{Xu}, any star-product in $C(\M)((\h))$
is equivalent to the one obtained via Fedosov's
procedure\footnote{See paper \cite{D}, in which this result was extended to
any smooth affine algebraic symplectic variety.}. In particular,
the cohomology class of the Fedosov-Weyl curvature (\ref{Weyl-c})
is a well-defined characteristic class of a star-product
in $C(\M)((\h))$. This characteristic class is referred to as
{\it the Deligne-Fedosov class}.

Let us prove here an important technical lemma
which might as well have an independent interest
\begin{lem}
\label{Tech}
Let $V$ be a vector bundle over the symplectic
manifold $\M$ and $D$ be a Fedosov connection
(\ref{F-conn}) with the form $r = A + dx^i \om_{ij}(x)y^j$
satisfying the normalization
condition (\ref{norm}).
If $q\in \EV$ is an endomorphism of $V$ then
$$
\pa^V q =0 \quad \Rightarrow \quad D q=0\,.
$$
In other words, if $q$ is $\pa^V$-flat then
the isomorphism $\la$ (\ref{la}) sends the
$q$ to itself.
\end{lem}
{\bf Proof.}
Since $q$ does not depend on $y's$ and
$\pa^V q=0$ we have that $\n q=\de q=0\,.$
Hence, $D q= \h^{-1} [r,q]$
and it suffices to show that the commutator
$[r,q]$ is vanishing. Since $D$
is nilpotent $D [r,q] = 0$.

The operator $\de^{-1}$ (\ref{del-1}) is not a derivation
of algebra $\Om^{\bul}(W(\EV))$. However, since $q$ does
not depend on $y's$ and $\de^{-1}r=0$,
$$
\de^{-1} [r,q]=0\,.
$$
Thus if $\eta =[r,q]\in \Om^1(W(\EV))$ we have
$$
\de \eta = \n \eta + \frac1{\h}[r,\eta]\,, \qquad
\de^{-1} \eta=0\,.
$$
Therefore, applying (\ref{Hodge}) to $\eta$ we get
$$
\eta = \de^{-1} (\n \eta + \frac1{\h}[r,\eta])\,.
$$
The latter equation has the unique vanishing solution
since $\de^{-1}$ raises the degree in $y$\,. Thus the lemma
is proven. $\Box$

\subsection{The Feigin-Felder-Shoikhet trace density map.}
Let as above $V$ be a vector bundle of rank $N$
over the symplectic manifold $\M$ and $*$ be a star-product
(\ref{star}) in the algebra $\EV((\h))$\,. Then,
\begin{defi} A trace density map $trd$ is a
$\bbC((\h))$-linear map
$$
trd\,:\, \EV((\h)) \mapsto H^{2n}(\M)((\h))
$$
vanishing on commutators
$$
trd(a * b - b * a)=0\,,
\qquad a,b \in \EV((\h))\,.
$$
\end{defi}
We will show that if the star-product $*$ in
$\EV((\h))$ is obtained via the Fedosov procedure
\cite{F1} then the Feigin-Felder-Shoikhet
cocycle (\ref{Te}) provides us with a natural
trace density map.

Let as above $D$ be the Fedosov connection (\ref{F-conn})
and $\la$ be the isomorphism (\ref{la}).
Since $\Te_{2n}^N$ (\ref{Te}) is a cocycle of $\mg= \mgl_N(\cW)$
relative to $\mh = \mgl_N\oplus \msp_{2n}$\,, we have the
following well-defined map
\begin{equation}
\label{Psi}
\Psi_{D}\,:\, \EV((\h)) \mapsto \Om^{2n}(\M)((\h))\,,
\qquad
\Psi_{D}(a)= \frac1{\h^n}\Te^N_{2n} (A,\dots, A,\la(a))\,,
\end{equation}
where $A\in \Om^1(W(\EV))$ is the one-form entering the definition
of the Fedosov connection $D$ (\ref{F-conn})\,.

We assemble the required properties of the map $\Psi_D$
in the following
\begin{teo}[Feigin-Felder-Shoikhet, \cite{FFS}]
\label{main} With the above notations the following
statements hold:\\
{\bf i)} If $V_1$ and $V_2$ are
two vector bundles over $\M$ and $D_1,~ D_2$ are
Fedosov connections on $W(End_{V_1})$ and
$W(End_{V_2})$, respectively, then for any $a\in End_{V_1}$
$$
\Psi_{D_1\oplus D_2} (a \oplus 0) = \Psi_{D_1}(a)\,,
$$
where $0$ stands for the trivial endomorphism
of $V_2$\,.\\
{\bf ii)} For any pair of endomorphisms $a,b \in \EV((\h))$
$$
\Psi_D(a * b - b * a)\in d \Om^{2n-1}(\M)((\h))\,,
$$
where $d$ is the De Rham differential.\\
{\bf iii)} Let $D$ and $\tD$ be two equivalent Fedosov
connections and $U\in I_N \oplus \G(W^1(\EV))$ be
the element establishing their equivalence in the
sense of (\ref{equiv1}). Then for any $a\in \EV((\h))$
$$
\Psi_{\tD}(\si (U^{-1}\circ \la(a)\circ U))
- \Psi_D (a)
\in d \Om^{2n-1}(\M)((\h))\,,
$$
where $\si$ is the projection (\ref{sigma}) and
$\la$ is the isomorphism (\ref{la}).\\
{\bf iv)} For the identity endomorphism $I_N\in \EV$
$$
\Psi_D(I_N) - \frac1{\h^n}Q_n(C(A,A), \dots, C(A,A))
\in d \Om^{2n-1}(\M)((\h))\,,
$$
where $Q_n\in ((S^n\mh)^*)^{\mh}$ is the $n$-th component
of the adjoint invariant form
(\ref{Qn}), $A\in \Om^1(W(\EV))$ is the one-form
entering the definition of the Fedosov
connection $D$ (\ref{F-conn})
and $C$ is the fiberwise ``curvature''
(\ref{curvature}) of the projection (\ref{pr})\,.
\end{teo}
{\bf Proof.} Statement $i)$ is obvious from the
construction of the cocycle $\Te^N_{2n}$ (\ref{phi}), (\ref{tau}),
(\ref{Te}). Statements $ii)$, $iii)$, and $iv)$
are proven in \cite{FFS} for the case when the
bundle $V$ is trivial. Statements $ii)$
and $iii)$ follow from the fact $\Te^N_{2n}$ is
a cocycle relative to $\mh = \mgl_N\oplus \msp_{2n}$
and statement $iv)$ is a consequence of (\ref{FFS})\,.
These results of \cite{FFS} are generalized to the
non-trivial bundle $V$ in a straightforward manner
since we always deal with the Lie algebra cochains
relative to $\mh = \mgl_N\oplus \msp_{2n}$\,. $\Box$

As an immediate consequence of this theorem we get that
\begin{cor} If $*$ is the star-product in $\EV((\h))$
obtained via the isomorphism $\la$ (\ref{la}) then
\begin{equation}
\label{trd}
trd(a) = [\Psi_D(a)] \,:\,
\EV((\h)) \mapsto H^{2n}(\M)((\h))
\end{equation}
is a trace density map.
\end{cor}
In what follows we refer to (\ref{trd}) as
the Feigin-Felder-Shoikhet (FFS) trace density map.

\section{The local version of the algebraic index theorem}
Let $\M$ be either a smooth real symplectic ma\-ni\-fold
of dimension $2n$ or a smooth affine algebraic variety (over $\bbC$)
of the complex dimension $2n$ endowed with an algebraic
symplectic form.  Let $*$ be a star-product in
the vector space $C(\M)((\h))$ of smooth (resp. regular) functions
on $\M$. Throughout this section we denote the algebra $(C(\M)((\h)),*)$
by $\bbA$ and the subalgebra $(C(\M)[[\h]],*)$ of $\bbA$ by $\bbA^+$

Due to remark $2$ after theorem \ref{F},
we may safely assume that $*$ is obtained by Fedosov's
procedure. Since in this
case the vector bundle is trivial the Fedosov star-product
depends only on the pair $(\pa^s, \Om_{\h}) $, where
$\pa^s$ is the symplectic connection on $\M$ and
$\Om_{\h}\in \h\Om^2(\M)[[\h]]$ is a series of closed
two-forms. Let $D= \n + \h^{-1}[A,\cdot\,] $ be the
corresponding Fedosov connection.

For any idempotent $P$ in the matrix algebra
$\mgl_N(\bbA^+)$ we assign the top degree
De Rham cohomology class
\begin{equation}
\label{class}
cl(P) = [\Psi_D(P)] \in H^{2n}(\M)((\h))\,,
\end{equation}
where $D$ is naturally extended to the
Fedosov connection on the Weyl algebra
bundle $W(End(L_N))$ associated with the trivial
bundle $L_N$ of rank $N$\,.

Due to statements $i)$ and $ii)$ (\ref{class})
gives a well-defined map
$$
cl \,:\, K_0(\bbA^+) \mapsto H^{2n}(\M)((\h))
$$
from the $K_0$-group of the algebra $\bbA^+$
to $H^{2n}(\M)((\h))$.
\begin{defi} The zeroth term $q=P|_{\h=0}$
of an idempotent $P\in \mgl_N(\bbA^{+})$ is
called the principal part of $P$\,.
\end{defi}
It is obvious that $q$ is an idempotent in the
matrix algebra $\mgl_N(C(\M))$ and the operation
of taking the principal part gives
a well-defined principal symbol map
\begin{equation}
\label{symb}
\Xi\,:\, K_0(\bbA^+) \mapsto
K_0(C(\M))\,.
\end{equation}

Due to the observation of L. Boutet de Monvel and
V. Guillemin \cite{BM-G} for any idempotent
$q$ in the matrix algebra $\mgl_N(C(\M))$ there exists
an idempotent $P\in \mgl_N(\bbA^+)$ whose principal part
is $q$. (See the explicit formula for $P$ in
\cite{F1}, eq. ($6.1.4$))

Using standard arguments of the index theory
we can prove that the cohomology class (\ref{class})
of an idempotent $P\in\mgl_N(\bbA^+)$ depends only
on its principal part. For our purposes we
need a slightly more general statement:
\begin{pred}
\label{only}
Let $V$ be a vector bundle over $\M$ and
$*$ be the Fedosov star-product in
$\EV[[\h]]$. Then if two
idempotents $P_1, P_2\in (\EV[[\h]], *)$
have the same principal part
$$
P_1\Big|_{\h=0} = P_2\Big|_{\h=0}\,,
$$
the cohomology classes
$[\Psi_D(P_1)]$ and $[\Psi_D(P_2)]$
coincide.
\end{pred}
{\bf Proof} is a straightforward generalization
of the arguments in the proof of theorem $6.1.3$ in
\cite{F1}\,. $\Box$

The precise dependence of the cohomology class
$cl(P)$ on the principal part $q$ of $P$ is given
by the local algebraic index theorem
\begin{teo}\label{index}
For any element $\Pi \in K_0(\bbA^+)$ the
cohomology class $cl(\Pi)$ coincides with the
top component of the cup product
$$
cl(\Pi) = \Big[ \,\hA(\M) \exp\left(-\frac{F}{\h}\right)\, ch(\Xi(\Pi))\,
\Big]_{2n}\,,
$$
of the $\hA$-genus of $\M$, the
exponent $\displaystyle e^{-\frac{F}{\h}}$ of the Deligne-Fedosov class
$$
F= \Big[ -\om + \Om_{\h} \Big]
$$
of the star-product in $\bbA^+$
and the Chern character  $ch(\Xi(\Pi))$
of the principal symbol $\Xi(\Pi)$ of $\Pi$\,.
\end{teo}
{\bf Proof.} Let $P$ be an idempotent in the matrix
algebra $\mgl_N(\bbA^+)$ and $q\in \mgl_N(C(\M))$ be the
principal part of $P$. As an idempotent endomorphism
of the trivial bundle $L_N$ or rank $N$, $q$ defines
a subbundle $E= Im\,q\subset  L_N$.

It suffices to prove that the
$2n$-form $\Psi_D(P)$ has the same cohomology
class as the $2n$-th component of the form
$$
\det\Big(\frac{R/2}{sinh\,(R/2)}
\Big)^{1/2}\exp\left( \frac{\om - \Om_{\h}}{\h}\right)
\, tr \exp(R^E)\,,
$$
where $R$ is the curvature form of the
symplectic connection $\pa^s$, $R^E$ is the
curvature form of the vector bundle $E$\,,
and $tr$ stands for the ordinary matrix
trace.

We start with an observation that $q$ is a flat section
of $End(L_{N})$ with respect to the connection
\begin{equation}
\label{pa-L}
\pa^{L_N} = d + [q(dq)-(dq)q,\cdot\,]\,.
\end{equation}

Let $D^{L_N}$ be the Fedosov connection
(\ref{F-conn}) corresponding
 to the initial symplectic connection
$\pa^s$ on $\M$, the initial series $\Om_h$ of
closed two-forms and the connection (\ref{pa-L})
on the trivial bundle $L_N$\,. By claim $3$ of theorem
\ref{F} there exists an element
$U\in I_N \oplus W^1(\EL)$ which establishes
the equivalence between the connections $D$ and $D^{L_N}$
in the sense of definition \ref{equiv}.
Therefore, by statement $iii)$ of theorem \ref{main}
$$
\Psi_{D}(P) - \Psi_{D^{L_N}}(\tP) \in
d\Om^{2n-1}(\M)((\h))\,,
$$
where $\tP= \si (U^{-1}\circ P \circ U)$ is
an idempotent in the algebra $(\EL, \tilde{*})$,
with the star-product $\tstar$ corresponding
to the Fedosov connection $D^{L_N}$\,.

Notice that, $\tP$ and $P$ have the
same principal part $\tP\Big|_{\h=0}=P\Big|_{\h=0}=q$.
Let us consider $q$ as an element in the algebra
$(\EL, \tilde{*})$\,. By claim $1$ of theorem
\ref{F} we may assume that $D^{L_N}$ satisfies
normalization condition (\ref{norm}). Therefore, by
lemma \ref{Tech} the element $q$ is flat with
respect to $D^{L_N}$, and hence $q$
is an idempotent of the algebra $(\EL, \tstar)$\,.
By proposition \ref{only} we may safely assume that
$\tP=q$\,.

We recall that $E= Im\,q\subset L_N$ is the subbundle
of $L_N$\,, corresponding to the idempotent $q$ and
denote by $\tE$ the subbundle $\tE\subset L_N$ corresponding to the
complementary idempotent $I_N-q\in \EL$\,. Thus, our trivial
vector bundle $L_N$ is isomorphic to the direct
sum $E\oplus \tE$\,. Let us fix an isomorphism $u$
$$
u\,:\, L ~\widetilde{\rightarrow}~ E\oplus \tE\,.
$$

Since $\tq= u q u^{-1}\in End_{E\oplus \tE}$
is the projector onto $E$ along $\tE$
\begin{equation}
\label{E-m}
\tq\, \Big|_{\G(E)((\h))} = I_m\,, \qquad
\tq\, \Big|_{\G(\tE)((\h))} = 0\,,
\end{equation}
where $m$ is the rank of $E$ and
$I_m$ is the identity endomorphism of $E$\,.

Due to remark $1$ after theorem \ref{F}
\begin{equation}
\label{conj}
D^{E\oplus \tE} = D^{L_N} + [u D^{L_N}(u^{-1}), \cdot\,]
\end{equation}
is the Fedosov connection on the Weyl algebra bundle
$W(End_{E\oplus \tE})$.

Furthermore, $\tq= u q u^{-1}$ is flat with respect to
$D^{E\oplus \tE}$, and hence the connection $D^{E\oplus \tE}$
preserves the following subspaces of
$End_{E\oplus \tE}((\h))$
$$
\{a\in End_{E\oplus \tE}((\h))~|~ a\,\Big|_{\G(E)((\h))}
\subset \G(E)((\h)),~ \ a\,\Big|_{\G(\tE)((\h))}=0 \}\,,
$$
$$
\{a\in End_{E\oplus \tE}((\h))~|~ a\,\Big|_{\G(\tE)((\h))}
\subset \G(\tE)((\h)),~ \ a\,\Big|_{\G(E)((\h))}=0 \}\,.
$$
The latter implies that $D^{E\oplus \tE}$
is a direct sum $D^E \oplus D^{\tE}$
of the Fedosov connections
\begin{equation}
\label{D-E}
D^E = \n^E  + \frac1{\h}[A^E, \cdot\,]\,,
\qquad
D^{\tE} = \n^{\tE}  + \frac1{\h}[A^{\tE}, \cdot\,]
\end{equation}
on the bundles $W(End_{E})$ and $W(End_{\tE})$, respectively.
Thus, due to (\ref{E-m}) and statement $i)$
of theorem \ref{main} we have
$$
\Psi_{D^{L_N}}(\tP) = \Psi_{D^E}(I_m)\,,
$$
where $I_m$ is as above the identity endomorphism of
the vector bundle $E$\,.

Using statement $iv)$ of theorem \ref{main} we get
that
\begin{equation}
\label{vot-ono}
\Psi_{D^E}(I_m) - \frac1{\h^n}\,Q_n(C(A^E,A^E), \dots, C(A^E,A^E))
\in d \Om^{2n-1}(\M)((\h))\,,
\end{equation}
where $Q_n\in ((S^n\mh)^*)^{\mh}$ is the $n$-th component
of the adjoint invariant form
(\ref{Qn}), $\mh= \mgl_N \oplus \msp_{2n}$\,,
$A^E\in \Om^1(W(End_{E}))$ is the one-form
of the Fedosov connection $D^E$ (\ref{D-E}),
and $C$ is the fiberwise ``curvature''
(\ref{curvature}) of the projection (\ref{pr}).

To compute $C(A^E, A^E)$ we mention that since
the Fedosov connection $D^{E\oplus\tE}$ is obtained
from the $D^{L_N}$ via conjugation (\ref{conj})
the connection $D^{E\oplus\tE}$ has the same
Fedosov-Weyl curvature (\ref{Weyl-c})
$$
\h(R+R^{E\oplus\tE}) +2 \n^{E\oplus\tE}
A^{E\oplus\tE} + \frac1{\h}
[A^{E\oplus\tE},  A^{E\oplus\tE}] =
-\om + \Om_{\h}\,.
$$
Here $A^{E\oplus\tE}$ is the one-form of the
Fedosov connection $D^{E\oplus\tE}$,
$\n^{E\oplus\tE}$ is the operator (\ref{nabla})
corresponding to the symplectic
connection $\pa^s$ and the connection $\pa^{E\oplus\tE}$
on $E\oplus \tE$ obtained from $\pa^{L_N}$ via the
isomorphism $u$\,. Finally, $R^{E\oplus\tE}$ is the
curvature form of the connection $\pa^{E\oplus\tE}$\,.

On the other hand the Fedosov connection $D^{E\oplus\tE}$
is a direct sum of the Fedosov connections $D^E$ and
$D^{\tE}$\,. Hence,
\begin{equation}
\label{the-same}
\h(R+ R^{E}) +2 \n^{E}
A^{E} + \frac1{\h}
[A^{E}, A^{E}] =
-\om + \Om_{\h}\,,
\end{equation}
where $R^E$ is the curvature form of the connection
$\pa^E$ on $E$ and $\n^E$ is the operator (\ref{nabla})
corresponding to the symplectic connection $\pa^s$ on $\M$
and the connection $\pa^E$\,.

Notice that, the projection $pr(A^E)$ may be non-zero.
However, using the trick of proposition \ref{absor}
we can absorb the terms that contribute to $pr(A^E)$
into the operator $\n^E$. In this way we change
both\footnote{Notice that,
both the new Christoffel form $\G$ and the new connection
form $\G^E$ become formal power series in $\h$\,.}
the symplectic connection $\pa^s$ on $\M$ and
the connection $\pa^E$ on $E$\,, and therefore we
change curvature forms $R$ and $R^E$. It is not hard to show that
after this rearrangement equation (\ref{the-same})
still holds.

Since $pr(A^E)$ is now vanishing $pr \n (A^E)= \n pr (A^E)=0$\,.
Thus, applying $pr$ to both sides of (\ref{the-same})
we get
\begin{equation}
\label{that}
C(A^E, A^E) = \h^2 R + \h^2 R_E - \h
(-\om + \Om_{\h})\,.
\end{equation}
Substituting (\ref{that}) into equation (\ref{vot-ono}) we complete the
proof of the theorem. $\Box$

\section{Concluding remarks.}
We would like to mention that our result (theorem \ref{index})
is derived from the purely algebraic fact (\ref{FFS})
about the cohomology of Lie algebra of matrices over
the Weyl algebra. In this respect our approach is very reminiscent
of finding ``local'' proofs of the Riemann-Roch-Hirzebruch
theorem \cite{ACKP}, \cite{BS}, \cite{BNT}, \cite{BNT1}, and \cite{FT}.
Thus, in \cite{ACKP} and \cite{BS} the Riemann-Roch theorem
for families of Riemann curves is deduced from
purely algebraic facts about the cohomology of the
Lie algebra of vector fields. In paper \cite{FT} this
idea is generalized to higher dimensions.

The Bressler-Nest-Tsygan theorem \cite{BNT1} or
the Riemann-Roch theorem for deformation quantizations, which we
already mentioned in the introduction, is also a local
statement based on purely algebraic relations between
Hochschild, cyclic, and Lie algebra (co)homology.
It seems that this result is the most general Riemann-Roch type
theorem which can be proven in the symplectic setting using the methods
of \cite{BNT1}.

A very similar statement to our result is proposed in recent
paper \cite{FLS} by B. Feigin, A. Losev, and B. Shoikhet.
In this paper the authors consider the algebra $\Dif(\cE)$ of
holomorphic differential operators acting on sections
of a holomorphic vector bundle $\cE$ over a compact
complex manifold $X$ of the complex dimension $n$.
They give a tractable notion
of Hochschild homology $HH_{\bul}$ of $\Dif(\cE)$
for which
$
\displaystyle HH_{\bul}(\Dif(\cE))= H^{2n-\bul}(X, \bbC)
$
and construct a map
\begin{equation}
\label{La}
\La\,:\, \Dif(\cE) \mapsto H^{2n}(X, \bbC)
\end{equation}
which induces an isomorphism $HH_0(\Dif(\cE)) \cong H^{2n}(X,\bbC)$\,.
Then, using the original Riemann-Roch-Hirzebruch
theorem and assuming that the Euler characteristic of $\cE$
is nonzero the authors of \cite{FLS} deduce that for any
$\cD\in \Dif(\cE)$, the pairing of $\La(\cD)$ with
the fundamental class $[X]$ coincides with the
super-trace $str(\cD)$ of $\cD$\,.

We would like mention that
the technique used in \cite{FLS} to define
the map (\ref{La}) is very similar to the
procedure of constructing the Hochschild cocycle
(\ref{tau})\,. This technique was
originally proposed in paper \cite{L} and we suspect that
it has an independent interest.

Notice that, since our proof of theorem \ref{index} is purely algebraic, it
can be generalized in a straightforward manner to the
setting of the symplectic Lie algebroids \cite{NT-Lie}\,.
In this case the local algebraic index theorem
relates the cohomology classes of the De Rham
complex associated with the corresponding
Lie algebroid.


\begin{thebibliography}{99}

\bibitem{ACKP} E. Arbarello, C. de Concini, V. Kac, and
C. Procesi, Moduli spaces of curves and representation
theory, Commun. Math. Phys. {\bf 117} (1988) 1-36.

\bibitem{AS} M. Atiyah and I. Singer, The index of elliptic
operators, I, III, Ann. of Math. {\bf 87}, 2 (1968)
484-530, 546-609.

\bibitem{BS} A. Beilinson and V. Schechtman, Determinant bundle
and Virasoro algebras, Commun. Math. Phys. {\bf 118}
(1988) 651-701.

\bibitem{BM-G} L. Boutet de Monvel and V. Guillemin,
The spectral theory of Toeplitz operators, Ann. Math. Stud.,
Vol. 99, Princeton University Press, 1981.

\bibitem{BNT} P. Bressler, R. Nest, and B. Tsygan,
Riemann-Roch theorems via deformation quantization.
I, Adv. Math. {\bf 167}, 1 (2002) 1-25.

\bibitem{BNT1} P. Bressler, R. Nest, and B. Tsygan,
Riemann-Roch theorems via deformation quantization.
II, Adv. Math. {\bf 167}, 1 (2002) 26-73.

\bibitem{CFS} A. Connes, M. Flato, and D. Sternheimer,
Closed star-products and cyclic cohomology, Lett. Math. Phys.
{\bf 24} (1992) 1-12.

\bibitem{D} V.A. Dolgushev, Hochschild cohomology versus De Rham
cohomology without formality theorems,  math.QA/0405177.

\bibitem{F1} B. Fedosov, Deformation quantization and index theory,
Akademie Verlag, Berlin, 1996.

\bibitem{FFS} B. Feigin, G. Felder, and B. Shoikhet,
Hochschild cohomology of the Weyl algebra and trace in deformation quantization,
math.QA/0311303.

\bibitem{FLS} B. Feigin, A. Losev, and B. Shoikhet,
Riemann-Roch-Hirzebruch theorem and to\-po\-lo\-gi\-cal
quantum mechanics, math.QA/0401400.


\bibitem{FT} B. Feigin and B. Tsygan, Riemann-Roch theorem and Lie algebra
cohomology I, Proceeding of the winter school on geometry and physics,
Srni, 1988, Suppl. Rend. Circ. Mat. Palermo, Ser, II, 21 (1989, 15-51)

\bibitem{K} M. Kontsevich, Deformation quantization of Poisson
manifolds, Lett. Math. Phys. {\bf 66} (2003) 157-216;
q-alg/9709040.

\bibitem{L} V. Lysov, Anticommutativity equation in topological
quantum mechanics, JETP Lett. 76 (2002) 724-727; hep-th/0212005.

\bibitem{RM} R. Melrose, Star products and local line
bundles, preprint, 2004, http://www-math.mit.edu/\~{}rbm/

\bibitem{NT} R. Nest and B Tsygan, Algebraic index theorem,
 Commun. Math. Phys. {\bf 172}, 2 (1995) 223-262.

\bibitem{NT-Lie} R. Nest and B Tsygan,
Deformations of symplectic Lie algebroids,
deformations of holomorphic symplectic structures, and index
theorems, Asian J. Math. {\bf 5}, 4 (2001), 599--635;
math.QA/9906020.


\bibitem{NT1} R. Nest and B Tsygan, Formal versus analytic index
theorems, IMRN, no. 11 (1996) 557-564.

\bibitem{Xu} P. Xu, Fedosov *-products and quantum momentum maps,
Commun. Math. Phys. {\bf 197}, 1 (1998) 167--197.
\end{thebibliography}
\end{document}